\documentclass[11pt,a4paper]{article}
\usepackage{amsmath,amssymb}

\usepackage{amsmath}

\usepackage{color}
\usepackage[colorlinks,linkcolor=blue,citecolor=red]{hyperref}

\DeclareMathSymbol{\bbbr}{\mathalpha}{AMSb}{"52}
\DeclareMathSymbol{\bbbc}{\mathalpha}{AMSb}{"52}

\newtheorem{theorem}{Theorem}

\newtheorem{conjecture}[theorem]{Conjecture}

\newtheorem{proposition}[theorem]{Proposition}

\begin{document}

\title{Quasilinear systems with linearizable characteristic webs}

\author{S.I.  Agafonov$^1$, E.V. Ferapontov$^2$ and V.S. Novikov$^2$}
     \date{}
     \maketitle
     \vspace{-5mm}
\begin{center}
$^1$ Departamento de Matem\'atica\\
UNESP-Universidade Estadual Paulista\\ S\~ao Jos\'e do Rio Preto\\ Brazil\\
  \ \\
$^2$Department of Mathematical Sciences \\ Loughborough University \\
Loughborough, Leicestershire LE11 3TU \\ United Kingdom \\
\ \\
e-mails: \\[1ex]  \texttt{agafonov@ibilce.unesp.br} \\
\texttt{E.V.Ferapontov@lboro.ac.uk}\\
\texttt{V.S.Novikov@lboro.ac.uk}\\
\end{center}

\vspace{1cm}

\begin{abstract}
We classify quasilinear systems in Riemann invariants  whose characteristic webs are  linearizable on every solution. Although the linearizability of an individual web is a rather nontrivial differential constraint, the requirement of linearizability of characteristic webs on {\it all} solutions imposes simple second-order constraints for the characteristic speeds of the system.
It is demonstrated that every such system with $n>3$ components can be transformed  by a  reciprocal transformation to $n$ uncoupled  Hopf equations.  All our considerations are local.

\bigskip


\bigskip

\end{abstract}

\newpage



\tableofcontents

\section{Introduction}

In this paper we investigate the geometry of characteristics of quasilinear systems in Riemann invariants
\begin{equation}
R^i_t=\lambda^i(R) R^i_x,
\label{n}
\end{equation}
$i=1, \dots, n$ (no summation over $i$). Systems of this form govern a wide range of problems in pure and applied mathematics.  Let us recall the basic  concepts needed to state our main results.

\medskip
\noindent {\bf Semi-Hamiltonian property}. System (\ref{n}) is called {\it semi-Hamiltonian} if its characteristic speeds $\lambda^i$ satisfy the   constraints
$$
\partial_k \left(\frac{\partial_j\lambda
^i}{\lambda^j-\lambda^i}\right)=\partial_j \left(\frac{\partial_k\lambda
^i}{\lambda^k-\lambda^i}\right),
$$
here $\partial_i=\partial/\partial {R^i}$. It was shown by Tsarev \cite{Tsar}  that this property is equivalent to integrability:  semi-Hamiltonian systems (\ref{n}) possess infinitely many conservation laws and commuting flows, and can be solved by the generalised hodograph method.

\medskip
\noindent {\bf Linear degeneracy}. System (\ref{n}) is said to be {\it linearly degenerate} if its characteristic speeds satisfy the conditions
$$
\partial_i\lambda^i=0,
$$
no summation,  $i=1, \dots, n$. Linear degeneracy is known to prevent the breakdown of smooth initial data, which is typical for genuinely nonlinear systems of type (\ref{n}),  and has been thoroughly investigated in the literature, see e.g. \cite{Roz, Liu, Serre}.

\medskip
\noindent {\bf Reciprocal transformations}. There exists a natural class of {\it reciprocal transformations} acting on systems of type (\ref{n}). These  are  non-local changes of the independent variables, $(x, t) \to (\tilde x, \tilde t)$, defined as
\begin{equation}\label{reciprocal}
 d\tilde x= Adt+Bdx, ~~~ d \tilde t=Mdt+Ndx,
 \end{equation}
 where the right-hand sides are two conservation laws of system (\ref{n}), that is, two 1-forms that are closed on every solution ($A, B, M, N$ are functions of $R$'s).
 The transformed system reads
 $$
 R^i_{\tilde t}=\tilde \lambda^i(R) R^i_{\tilde x},
 $$
 where the transformed characteristic speeds are given by
 \begin{equation}\label{newspeed}
\tilde \lambda^i=\frac{\lambda^iB-A}{M-\lambda^iN}.
 \end{equation}
Reciprocal transformations are known to preserve both the semi-Hamiltonian property and the linear degeneracy \cite{Fer1}.

\medskip
\noindent {\bf Characteristics}. {\it Characteristic curves} of the $i$-th family are defined by the equation $dx+\lambda^idt=0$. Altogether, characteristics form an $n$-web (that is, $n$ one-parameter families of curves) on every solution. We refer to \cite{BB, Bla} for an introduction to the web geometry.

\medskip
\noindent {\bf Parallelizable webs}. An $n$-web is said to be {\it parallelizable}   if it  is {\it locally diffeomorphic}  to $n$ families of parallel lines  (parallelizable $3$-webs are also known as {\it hexagonal}).
The following result provides a link between the above concepts:

 \begin{theorem} \cite{Fer2, Fer3} The following conditions are equivalent:

\noindent (a) System (\ref{n}) has a parallelizable characteristic web on every solution.

\noindent (b) System  (\ref{n}) can be linearized  by a reciprocal transformation.

\noindent (c) System (\ref{n}) satisfies the following conditions:

\begin{itemize}
\item semi-Hamiltonian property;

\item linear degeneracy;

\item every quadruple of characteristic speeds has a  constant cross-ratio   (for $n=3$ this condition becomes redundant).

\end{itemize}
\end{theorem}

\noindent {\bf Linearizable webs}. An $n$-web is said to be {\it linearizable} (rectifiable)   if it  is {\it locally diffeomorphic}  to $n$ families of  lines, not necessarily parallel. We emphasize that the condition of linearizability is far more subtle than that of parallelizability, see \cite{Bk, Bg, BB,  Henaut, GMS, AGL, GL, Ag} and references therein for a discussion of the linearizability problem.

\medskip

The  aim of this paper is to establish an analogue of Theorem 1 for systems (\ref{n}) whose characteristic webs are linearizable on every solution. Our main observation is that, although the linearizability of an individual web is a  rather nontrivial  differential constraint, the requirement of linearizability of characteristic webs on {\it all} solutions leads to simple second-order differential constraints for the characteristic speeds, see Theorem \ref{rec4} below.
In particular, any such system is reciprocally related to $n$ uncoupled Hopf equations, $R^i_t=f^i(R^i)R^i_x$.

\section{Linearizability of a planar web}
Any projective transformation takes a linear web into a linear web. Thus, the linearizability of a web can be expressed in terms of projective differential invariants of the linearizing map, namely, its Schwarzian derivative \cite{Sp}. This approach can be traced back to the pioneering work of Bol \cite{Bk,Bg}. The following form of the linearizability criterion will be  most convenient for our purposes:
\begin{theorem}[H\'enaut,  \cite{Henaut}]
A planar $n$-web formed by integral curves of vector fields
$$\{ V_i=\partial_t-\lambda^i (t,x)\partial_x\},$$
 is linearizable if and only if there exists a solution $E,F,G,H$ of the following system of PDEs,
\begin{equation}
\begin{array}{c}\label{compatibility}
2G_{tx}+H_{xx}+F_{tt}-6GG_x+2HE_t+EH_t+3FH_x-3GF_t+3HF_x=0,\\
 \\
G_{xx}+E_{tt}+2F_{tx}+3EG_t-3FG_x+3GE_t+HE_x+2EH_x-6FF_t=0,
\end{array}
\end{equation}
subject to the constraints
 \begin{equation}\label{cubH}
E(\lambda^i)^3+3F(\lambda^i)^2+3G\lambda^i+H= V_i(\lambda^i), \ \ \ i=1,..., n.
 \end{equation}
\end{theorem}

\medskip

\noindent {\bf Remark.} The functions $E,F,G,H$ are linearly independent components of the 2-dimensional Schwarzian derivative of the linearizing map \cite{Ag}.
For $n\geq 4$, relations (\ref{cubH}) uniquely define $E,F,G,H$, and equations (\ref{compatibility}) lead to explicit  second-order constraints for $\lambda^i$. On the contrary, for $n=3$, relations  (\ref{cubH}) are only sufficient to determine, say,  $E, F, G$ in terms of $H$, so that relations  (\ref{compatibility}) give an over-determined second-order system for $H$. The analysis of this system is quite involved, in particular, differential constraints for $\lambda^i$ appear at differential order six and higher. This explains why we treat the cases $n\geq 4$ and $n=3$   separately.

\section{$n\times n$ systems with linearizable characteristics ($n\geq 4$)}

The main result of this paper is the following Theorem.

\begin{theorem}\label{rec4} For $n\geq 4$ the following conditions are equivalent:

\noindent (a) System (\ref{n}) has a linearizable characteristic web on every solution.

\noindent (b) System  (\ref{n}) can be transformed by a reciprocal transformation to  $n$ uncoupled Hopf  equations.

\noindent (c) Characteristic speeds of system (\ref{n})   satisfy the following conditions:
\begin{itemize}
\item For every quadruple of pairwise distinct indices $i, j, k, l$ one has a relation
\begin{equation}
a_{ij}(\lambda^k-\lambda^l)+a_{kj}(\lambda^l-\lambda^i)+a_{lj}(\lambda^i-\lambda^k)=0,
\label{1st}
\end{equation}
here  $a_{ij}=\frac{\partial_j\lambda^i}{\lambda^j-\lambda^i}$.
This allows one to introduce the parametrization $a_{ij}= p_j\lambda^i+q_j$.

\item  The 1-forms
\begin{equation}
\omega_{11}=\sum_{i=1}^n p_i\lambda^idR^i, ~~ \omega_{12}=\sum_{i=1}^n q_i\lambda^idR^i, ~~ \omega_{21}=\sum_{i=1}^n p_idR^i, ~~ \omega_{22}=\sum_{i=1}^n q_idR^i,
\label{omega}
\end{equation}
 satisfy the $gl(2)$ structure equations,
\begin{equation}\label{structure}
d\omega_{ab}=\sum^2_{c=1} \omega_{ac}\wedge \omega_{cb},
\end{equation}
$a, b=1, 2$, which are equivalent to $\partial_jp_i=a_{ij}p_i, \ \partial_jq_i=a_{ij}q_i.$
\end{itemize}

\noindent A system satisfying either of the equivalent conditions (a), (b) or (c) is automatically semi-Hamiltonian.

\end{theorem}

\centerline{\bf Proof:}

\medskip The equivalence of (b) and (c) can be seen as follows. Consider  reciprocal transformation (\ref{reciprocal}),
note the relations $\partial_iA=\lambda^i\partial_iB, \ \partial_iM=\lambda^i\partial_iN$.  Requiring that the transformed characteristic speed (\ref{newspeed})
depends on the variable $R^i$ only, namely $\partial_j\tilde \lambda^i=0$ for every $j\ne i$,  we obtain all first-order partial derivatives of $A, B, M, N$.
Comparing the relations $\partial_j\tilde \lambda^i=\partial_j\tilde \lambda^k=\partial_j\tilde \lambda^l=0$ we  obtain  first-order relations (\ref{1st}). This allows one to set $a_{ij}= p_j\lambda^i+q_j$, leading to
$$
dA=A\omega_{11}+B\omega_{12}, ~~~ dB=A\omega_{21}+B\omega_{22},
$$
$$
dM=M\omega_{11}+N\omega_{12}, ~~~ dN=M\omega_{21}+N\omega_{22},
$$
where $\omega_{ab}$ are as in (\ref{omega}).  The compatibility conditions of these relations are nothing but  $gl(2)$ structure equations (\ref{structure}). A direct calculation shows  that they are equivalent to $\partial_jp_i=a_{ij}p_i, \ \partial_jq_i=a_{ij}q_i.$
By construction, in the new independent variables $\tilde x, \ \tilde t$ the system reduces to $n$ uncoupled Hopf equations, $R^i_{\tilde t}=f^i(R^i)R^i_{\tilde x}$, where $f^i$ are arbitrary functions
(possibly, constants). Note that non-constant $f^i$ can  be reduced to $R^i$ via a reparametrisation of Riemann invariants.

\medskip

Since characteristics of a Hopf equation are straight lines, and reciprocal transformations are nothing but non-local changes of variables depending on a solution, this also establishes the implication (b) $\implies$ (a)
(equivalently, (c) $\implies$ (a)).

\medskip

As for the less elementary implication, (a) $\implies$ (c), let  us first consider the case $n=4$. Equations  (\ref{cubH}) give $E,F,G,H$ as functions of the characteristic speeds  $\lambda^i$ and their derivatives, $\lambda^i_x=\sum_k\lambda^i_kR^k_x$, $\lambda^i_t=\sum_k\lambda^i_k\lambda^kR^k_x$ (we fix a  solution so that $R^i$ and $\lambda^i$ become functions of $x,t$). Substituting these expressions into  (\ref{compatibility}) we obtain polynomial equations in the differential variables $R^i_x,R^i_{xx},R^i_{xxx}$. Since the characteristic 4-web is required to be linearizable on {\it every} solution, the equations split with respect to these variables. Thus, equating to zero  coefficients at $R^i_{xxx}$  we get  first-order relations (\ref{1st}). Further, taking coefficients at $R^i_{xx}R^j_x$ and differentiating relations (\ref{1st}) with respect to $R^k$ yields all second-order relations for $\lambda^i$ that are equivalent to  $gl(2)$ structure equations (\ref{structure}).

Now the general case $n>4$ readily follows. Let as fix four pairwise distinct indices $i,j,k,l$, say $1, 2, 3, 4$. Consider special solutions to system (\ref{n}) such that $R^s=const, \ s>4$. This  reduces system (\ref{n}) to a  4-component system for $R^1, \dots, R^4$, with linearizable characteristic 4-webs. Therefore we have all relations (\ref{1st}), as well as all other necessary conditions involving indices $1, 2, 3, 4$. The rest follows from the fact that every linearizability condition involves maximum $4$ distinct indices.

\medskip
 Finally, the semi-Hamiltonian property follows from the fact that a system of uncoupled Hopf equations is automatically semi-Hamiltonian, and reciprocal transformations preserve the semi-Hamiltonian property.

\hfill $\Box$

\bigskip

\noindent {\bf Remark.} Note that the equations for $A, B$ and $M, N$ uncouple and the forms $\omega_{ab}$ verify  $gl(2)$ structure equations. Therefore $(A, B)$ and $(M, N)$ can be interpreted as sections of a flat 2-dimensional vector bundle over the hodograph space.

\bigskip

\noindent {\bf Example.} The simplest example of an $n$-component system with linearizable characteristics is provided by $n$ uncoupled Euler equations,
$$
R^i_t=R^iR^i_x.
$$
Remarkably, it appears as the modulational Whitham system for the Benjamin-Ono equation \cite{DK}. Applying to this system  reciprocal transformations one obtains  {generic} systems with linearizable characteristics. Explicitly, let us take two conservation laws,
$$
\begin{array}{c}
Adt+Bdx=[\sum ((f^k)' R^k-f^k)]dt+[\sum (f^k)']dx, \\
\ \\
Mdt+Ndx=[\sum ((g^k)' R^k-g^k)]dt+[\sum (g^k)']dx,
\end{array}
$$
where $f^k(R^k)$ and $g^k(R^k)$ are arbitrary functions, and prime indicates derivative. The transformed characteristic speeds take the form
$$
\tilde \lambda^i=\frac{\lambda^iB-A}{M-\lambda^iN}=-\frac{\sum [(f^k)'(R^i-R^k)+ f^k]}{\sum [(g^k)'(R^i-R^k)+ g^k]}.
$$
This formula gives characteristic speeds of {\it generic} systems with linearizable characteristics. Examples of this type appeared as hydrodynamic reductions of integrable hydrodynamic chains in \cite{Pavlov}. Degenerations can be obtained by replacing some of the Euler equations  by  linear equations, $R^i_t=c^iR^i_x$,
$c^i=const$. In particular, starting with $n$ linear equations and applying reciprocal transformations one obtains
$$
\tilde \lambda^i=-\frac{\sum(c^i-c^k)f^k(R^k)}{\sum(c^i-c^k)g^k(R^k)},
$$
which is a general formula for characteristic speeds of linearly degenerate semi-Hamiltonian systems with constant cross-ratios  \cite{Fer4}. In this limiting case the characteristic web is parallelizable on every solution.

\section{$3\times 3$ systems with linearizable characteristics}

Even though for $n=3$ relations (\ref{1st}) are vacuous,  we can always  represent $a_{ij}$ in the form $a_{ij}= p_j\lambda^i+q_j$. Explicitly,
$$p_j= \frac{\textstyle a_{ij}-a_{kj}}{\textstyle \lambda^i-\lambda^k}, \ \ q_j= \frac{\textstyle a_{kj}\lambda^i-a_{ij}\lambda^k}{\textstyle \lambda^i-\lambda^k}, $$
$ i, k\ne j$. Let us introduce the forms $\omega_{ab}$ by formulae
(\ref{omega}), where the summation is from 1 to 3. We have
$$
\begin{array}{c}
\omega_{11}=\frac{\textstyle a_{31}-a_{21}}{\textstyle \lambda^3-\lambda^2}\lambda^1dR^1+
\frac{\textstyle a_{12}-a_{32}}{\textstyle \lambda^1-\lambda^3}\lambda^2dR^2+\frac{\textstyle a_{23}-a_{13}}{\textstyle \lambda^2-\lambda^1}\lambda^3dR^3, \\
\ \\
\omega_{12}=\frac{\textstyle{ a_{21}} \lambda^3-\textstyle{a_{31}}\lambda^2}{\textstyle \lambda^3-\lambda^2}\lambda^1 dR^1+
\frac{\textstyle{a_{32}}\lambda^1-\textstyle {a_{12}}\lambda^3}{\textstyle \lambda^1-\lambda^3}\lambda^2 dR^2+\frac{\textstyle{a_{13}}\lambda^2-\textstyle{ a_{23}}\lambda^1}{\textstyle \lambda^2-\lambda^1}\lambda^3 dR^3,\\
\ \\
\omega_{21}=\frac{\textstyle a_{31}-a_{21}}{\textstyle \lambda^3-\lambda^2}dR^1+
\frac{\textstyle a_{12}-a_{32}}{\textstyle \lambda^1-\lambda^3}dR^2+\frac{\textstyle a_{23}-a_{13}}{\textstyle \lambda^2-\lambda^1}dR^3, \\
\ \\
\omega_{22}=\frac{\textstyle {a_{21}}\lambda^3-\textstyle {a_{31}}\lambda^2}{\textstyle \lambda^3-\lambda^2}dR^1+
\frac{\textstyle \textstyle {a_{32}}\lambda^1-\textstyle {a_{12}}\lambda^3}{\textstyle \lambda^1-\lambda^3}dR^2+\frac{\textstyle {a_{13}}\lambda^2-\textstyle {a_{23}}\lambda^1}{\textstyle \lambda^2-\lambda^1}dR^3.
\end{array}
$$

\begin{proposition}\label{propositionn3}
System  (\ref{n}) can be transformed by a reciprocal transformation to 3 uncoupled Hopf equations if and only if
the forms $\omega_{ab}$ satisfy  $gl(2)$ structure equations (\ref{structure}).
\end{proposition}

\centerline{\bf Proof:}

\medskip  This proposition claim can be proved exactly as in Theorem (\ref{rec4}). Namely, requiring $\partial_j\tilde \lambda^i=0$ for any $j\ne i$, we obtain all first-order partial derivatives of $A, B, M, N$ in the form
$$
dA=A\omega_{11}+B\omega_{12}, ~~~ dB=A\omega_{21}+B\omega_{22},
$$
$$
dM=M\omega_{11}+N\omega_{12}, ~~~ dN=M\omega_{21}+N\omega_{22},
$$
where $\omega_{ab}$ are as above. The compatibility conditions of these relations are the $gl(2)$ structure equations.

\hfill $\Box$

Due to the complexity of linearizability conditions for 3-webs, we were unable to prove the  analogue of Theorem (\ref{rec4}) for $n=3$. Thus, we can only formulate the following conjecture.

\begin{conjecture}\label{conjecture} For  3-component system (\ref{n}), the following conditions are equivalent:

\noindent (a) System (\ref{n}) has a linearizable characteristic 3-web on every solution.

\noindent (b) System  (\ref{n}) can be transformed by a reciprocal transformation to 3 uncoupled Hopf equations.

\noindent (c) The forms $\omega_{ab}$ satisfy  $gl(2)$ structure equations (\ref{structure}).

\noindent Any system satisfying either of the above conditions is automatically semi-Hamiltonian.
\end{conjecture}

\centerline{\bf Discussion:}

\medskip

The equivalence of (b) and (c) is proved in Proposition \ref{propositionn3}.

\medskip

The implication (b) $\implies$ (a), equivalently, (c) $\implies$ (a), is also straightforward: characteristics of a Hopf equation are straight lines, and reciprocal transformations are  non-local changes of variables depending on a solution.

\medskip

The main problem is the converse implication, (a) $\implies$ (c), which is a highly nontrivial computational challenge due to the complexity of linearizability conditions for 3-webs. This calculation seems to be out of reach for the modern computer algebra systems.
\hfill $\Box$

\bigskip

\noindent {\bf Remark.} The Gronwall conjecture  \cite{Gn} states that, modulo projective equivalence,  a linearizable non-hexagonal 3-web has a unique linear representation. While this is still open at the level of individual webs (for partial results see \cite{Bg,Wg,Ag}), our results show that the natural analogue of this statement holds at the level of systems. Namely, if a system can be decoupled by a reciprocal transformation, then the transformation  is unique up to linear changes of $x$ and $t$.

\section{Concluding remarks}

It remains a considerable computational challenge to establish the implication $(a) \implies (c)$ of Conjecture \ref{conjecture}.
The difference with the case $n>3$ can be explained geometrically as follows. Any 4-subweb of a planar $n$-web  defines a unique projective connection$\nabla$. The 4-subweb is linearizable if and only if all web leaves are geodesics of $\nabla$, and the curvature of the connection vanishes. In fact, the functions $E,F,G,H$ defined by equations (\ref{cubH}) (where $i$ now runs over the indices of the 4-subweb) are, up to constant factors, the so-called Thomas coefficients of $\nabla$ (see \cite{Pl}). The flatness of $\nabla$ manifests itself in equations (\ref{compatibility}).

\section*{Acknowledgements}

This research was supported by FAPESP grant \#2014/17812-0. SIA thanks the Mathematical Department of
Loughborough University, where this study was initiated, for a kind hospitality.

\end{document}